\documentclass[journal]{IEEEtran}

\ifCLASSINFOpdf

\else

\fi
\usepackage[cmex10]{amsmath}

\usepackage{amssymb}

\usepackage{graphics,graphicx,epic,eepic,epsfig,times,units} 

\usepackage{float}

\usepackage{psfrag}

\usepackage{ifpdf}

\usepackage{cite}

\usepackage{siunitx}

\usepackage{multirow}

\usepackage{xcolor}

\usepackage{algorithm,algorithmic}
\usepackage{tabularx}

\begin{document}

\title{A Distributed Algorithm for Demand Response with Mixed-Integer Variables}

\author{Sleiman~Mhanna,~\IEEEmembership{Student~MIEEE,}
        Archie~C.~Chapman,~\IEEEmembership{~MIEEE}
        and~Gregor~Verbi\v{c},~\IEEEmembership{Senior~MIEEE,}
}


\maketitle

\begin{abstract}
This letter presents a distributed algorithm for aggregating a large number of households with mixed-integer variables and intricate couplings between devices. The proposed distributed gradient algorithm is applied to the double smoothed dual function of the adopted DR model. Numerical results show that, with minimal parameter adjustments, the convergence of the dual objective exhibits a very similar behavior irrespective system size.
\end{abstract}


\IEEEpeerreviewmaketitle

\section{Introduction}

\IEEEPARstart{E}{fficient} load scheduling and aggregation is a problem of growing importance in the area of demand response (DR). However, this problem is particularly difficult for two main reasons. First, many household electrical devices have discrete operating points that can only be represented by mixed-integer variables, and moreover, household device uses are often coupled, thus giving household electricity demand a combinatorial structure. Therefore, solving this problem centrally may spell
intractability when the number of households is large. Second, solving this problem centrally would require sending all of the households' private information to the aggregator, which entails substantial communication overhead and privacy concerns.

Against this background, this letter proposes a distributed gradient algorithm applied to a double smoothed dual function. This work is not the
first attempt to solve a DR problem with mixed-integer variables in a distributed fashion. This problem is decomposed in terms of devices and solved in a distributed fashion using the proximal bundle method in \cite{DRwithMIC}. In contrast to \cite{DRwithMIC}, the DR problem in this work is decomposed in terms of households. Doing so, allows for a more expressive household model, which can incorporate the intricate couplings between storage devices, appliances and renewable energy resources.

\section{DR model and problem description}\label{sec:DRmodel}

The model comprises a set of agents $\mathcal{I}:=\{0,1,2, \ldots ,I\}$, where $0$ is the aggregator and each $i \neq 0$ is a household agent. Let $x^{t}_{i} \in \mathbb{R_{+}}$ be the demand for electric energy for every agent $i \neq 0$ at time-slot $t\in \mathcal{T}:=\{\tau,\ldots,\tau+T-1\}$. The aggregator faces a set of cost functions $C^{t}:\mathbb{R}_{+} \mapsto \mathbb{R}_{+}$, where $C^{t}\left(x^{t}_{0}\right)$ is the cost of supplying $x^{t}_{0}$ units of energy to the households at time-slot $t$.

Given the households' feasible schedule sets $X_{i\neq 0}$ and their demand profile $\boldsymbol{x_i}=\left[x_{i}^{\tau}, \ldots , x_{i}^{\tau+T-1} \right]$, the aggregator can (centrally) minimise the total energy cost per scheduling horizon $T$ by solving the following problem:
\begin{subequations}\label{eq:mincost}
\begin{align}
 & \min_{\boldsymbol{x_{i}} \in X_{i}} && \sum_{ t\in \mathcal{T}} C^{t}\left(x^{t}_{0}\right), \\
 & \text{subject to} && \sum_{ i\in \mathcal{I} \setminus 0} x^{t}_{i}=x^{t}_{0}, \qquad t \in \mathcal{T}. \label{eq:coupling}
\end{align}
\end{subequations}
The local constraints of agents $i \neq 0$ arise from the operating modes of different flexible loads including \emph{interruptible} (e.g. PHEVs, pool pumps) and \emph{non-interruptible} (e.g. washing machines, dishwashers) loads (as in \cite{DRwithMIC,FaithfulMDinDR}). Problem \eqref{eq:mincost} can also be written as $\mathcal{P}^{*}=\inf_{\boldsymbol{x} \in X} \left\{C\left(\boldsymbol{x}\right):A_{c}\boldsymbol{x}=\boldsymbol{0}\right\},$ where $\boldsymbol{x}=\left\{\boldsymbol{x_i}\right\}_{i\in \mathcal{I}}$, $X=\prod_{ i\in \mathcal{I}} X_{i}$, and $A_{c}$ is the coupling constraint matrix.

Problem \eqref{eq:mincost} is a mixed-integer program (MIP) that belongs to the class of NP-hard problems which are notorious for tending to be intractable (if solved centrally) when they grow in size. However, relaxing the coupling constraints \eqref{eq:coupling}, through the Lagrangian relaxation method, bestows a separable structure on problem \eqref{eq:mincost}. The problem can then be decomposed into $I+1$ independent subproblems that can be solved in parallel. The partial Lagrangian of problem \eqref{eq:mincost} can be written as
\begin{align} 
	\mathcal{L}\left(\boldsymbol{x},\boldsymbol{\lambda} \right) & =\sum_{ t\in \mathcal{T}} C^{t}\left(x^{t}_{0}\right)+\sum_{ t\in \mathcal{T}} \lambda^{t} \left(\sum_{ i\in \mathcal{I} \setminus 0} x^{t}_{i} - x^{t}_{0}\right),
\end{align}
where $\boldsymbol{\lambda}=\left[\lambda^{\tau}, \ldots , \lambda^{\tau+T-1} \right]$ is the vector of Lagrange multipliers.
Therefore, the Lagrange dual function would be
\begin{align}
 &\hspace{-2 mm} \mathcal{D}\left(\boldsymbol{\lambda}\right) =\inf_{\boldsymbol{x_i} \in X_{i}} \mathcal{L}\left(\boldsymbol{x},\boldsymbol{\lambda} \right)=\mathcal{D}_{0}\left(\boldsymbol{\lambda}\right)+\sum_{ i\in \mathcal{I}\setminus 0} \mathcal{D}_{i}\left(\boldsymbol{\lambda}\right), \text{~where} \\
 &\hspace{-2 mm}\mathcal{D}_{0}\left(\boldsymbol{\lambda}\right)=\inf_{\boldsymbol{x_{0}} \in X_{0}} \sum_{ t\in \mathcal{T}} \left(C^{t}\left(x^{t}_{0}\right) - \lambda^{t}x^{t}_{0} \right), \text{~and} \\
 &\hspace{-2 mm}\mathcal{D}_{i}\left(\boldsymbol{\lambda}\right)=\inf_{\boldsymbol{x_i} \in X_{i}} \sum_{ t\in \mathcal{T}} \lambda^{t} x^{t}_{i} , \quad i \in \mathcal{I} \setminus 0. \label{eq:subproblems}
\end{align}
Finally, the dual problem is $\max_{\boldsymbol{\lambda} \succeq \boldsymbol{0}} \ \ \mathcal{D}\left(\boldsymbol{\lambda}\right)$.
However, in this DR scenario, the concave dual function $\mathcal{D}\left(\boldsymbol{\lambda}\right)$ is typically nondifferentiable. Specifically, as the subproblems in \eqref{eq:subproblems} can have multiple optimal solutions for a given vector $\boldsymbol{\lambda}$, the dual function $\mathcal{D}\left(\boldsymbol{\lambda}\right)$ can be nonsmooth. Consequently, applying a conventional gradient method \cite{subgradientmethods} to this problem would exhibit a very slow convergence.

\section{Double smoothing method}\label{sec:doublesmoothing}
One way to obtain a smooth approximation of $\mathcal{D}\left(\boldsymbol{\lambda}\right)$ is to modify the subproblems in \eqref{eq:subproblems} to ensure a unique optimal solution for every $\boldsymbol{\lambda}$. The dual function is modified as follows:
\begin{align}
	\mathcal{D}_{\mu} \left(\boldsymbol{\lambda}\right) &=\mathcal{D}_{0}\left(\boldsymbol{\lambda}\right)+\sum_{ i\in \mathcal{I} \setminus 0} \mathcal{D}_{i,\mu} \left(\boldsymbol{\lambda}\right),  \text{~where} \label{eq:smootheddual} \\ 
	\mathcal{D}_{i,\mu} \left(\boldsymbol{\lambda}\right) &=\inf_{\boldsymbol{x_i} \in X_{i}} \left(\sum_{ t\in \mathcal{T}} \lambda^{t} x^{t}_{i} + \frac{\mu}{2}\left\|\boldsymbol{x_i}\right\|^2\right), \ \ i \in \mathcal{I} \setminus 0, 
\end{align}
and $\mu > 0$ is a smoothness parameter.
The modified dual function $\mathcal{D}_{\mu}\left(\boldsymbol{\lambda}\right)$ is smooth and its gradient $\nabla \mathcal{D}_{\mu}\left(\boldsymbol{\lambda}\right)=A_{c} \boldsymbol{x}_{\mu,\boldsymbol{\lambda}}$, where $\boldsymbol{x}_{\mu,\boldsymbol{\lambda}}$ denotes the unique optimal solution of problem \eqref{eq:smootheddual}, is Lipschitz-continuous with Lipschitz constant $L_{\mu}=\frac{\left\|A_{c}\right\|^2}{\mu}$.

The aim of this smoothing is to obtain a Lipschitz-continuous gradient for which efficient smooth optimisation methods can be applied \cite{intoncvxopt}. However, despite having a good convergence rate $\mathcal{D}_{\mu}\left(\boldsymbol{\lambda}^{*}\right)-\mathcal{D}_{\mu}\left(\boldsymbol{\lambda}_{k}\right)$ at iteration $k$ when applying a fast gradient method, the same good rate of convergence does not apply to $\left\|\nabla \mathcal{D}_{\mu}\left(\boldsymbol{\lambda}_{k}\right)\right\|$. 

Since the aim is not only to efficiently solve the dual problem but also to recover a nearly feasible solution to the primal \cite{doublesmoothing}, a second smoothing is applied to the dual function to make it strongly concave. The new dual function is written as
\begin{align}
	\mathcal{D}_{\mu,\kappa} \left(\boldsymbol{\lambda}\right) =\mathcal{D}_{0}\left(\boldsymbol{\lambda}\right)+\sum_{ i\in \mathcal{I} \setminus 0} \mathcal{D}_{i,\mu} \left(\boldsymbol{\lambda}\right)-\frac{\kappa}{2}\left\|\boldsymbol{\lambda}\right\|^{2},
\end{align}
which is strongly concave with parameter $\kappa > 0$, and whose gradient $\nabla \mathcal{D}_{\mu,\kappa}\left(\boldsymbol{\lambda}\right) = A_c \boldsymbol{x}_{\mu,\boldsymbol{\lambda}} - \kappa \boldsymbol{\lambda}$ is Lipschitz-continuous with constant $L_{\mu,\kappa}=L_{\mu}+\kappa$. 

\section{Fast gradient algorithm}\label{sec:algorithm}
The fast gradient method involves two multiplier updates,	$\boldsymbol{\lambda}_{k+1}=\hat{\boldsymbol{\lambda}}_{k}+\frac{1}{L_{\mu^k,\kappa^k}^{k}} \nabla \mathcal{D}_{\mu^k,\kappa^k}\left(\hat{\boldsymbol{\lambda}}_{k}\right),$ and $\hat{\boldsymbol{\lambda}}_{k+1}=\boldsymbol{\lambda}_{k+1}+\beta^{k} \left(\boldsymbol{\lambda}_{k+1}-\boldsymbol{\lambda}_{k}\right),$ where $\beta^{k}=\left(\sqrt{L_{\mu^k,\kappa^k}^{k}}-\sqrt{\kappa^{k}}\right)\left(\sqrt{L_{\mu^k,\kappa^k}^{k}}+\sqrt{\kappa^{k}}\right)^{-1}$.

The parameters of the algorithm are set as follows, $\mu^{k+1}=\alpha^{k+1}/D_X, \text{ and } \kappa^{k+1}= \mathrm{e}^{\left(\log\left(\kappa^{\text{maxiter}}/ \kappa^{1}\right)/\text{maxiter}\right)} \kappa^{k}$, where $\alpha^{k+1}= \mathrm{e}^{\left(\log\left(\alpha^{\text{maxiter}}/\alpha^{1}\right)/\text{maxiter}\right)} \alpha^{k}$, and $\text{maxiter}$ is the maximum number of iterations. The distributed algorithm is described in Algorithm~\ref{algorithm}.

 \begin{algorithm}
 \caption{: Distributed algorithm}
  \algsetup{linenosize=\tiny}
  \scriptsize
 \begin{algorithmic}[1]
 \renewcommand{\algorithmicrequire}{\textbf{Parameters:}}
 \REQUIRE $\boldsymbol{\lambda}_{1} \succeq \boldsymbol{0} $, $\kappa^{1} > 0$, $\kappa^{\text{maxiter}}=0.0001$, $\hat{\mu}^{\text{min}} \in \left[0.0001,0.005\right]$, $\text{maxiter} \in \left\{500,1000\right\}$.


 \STATE \textbf{Initialisation:} Households compute $D_{X_i}=\min \left\{ \frac{1}{2} \left\| \boldsymbol{x_i} \right\|^{2} : \boldsymbol{x_i} \in X_{i} \right\},$ and send it to the aggregator which computes $D_X=\sum_{i \in \mathcal{I} \setminus 0} D_{X_i}$ and sets \\ $\mu^{1}=\alpha^{1}/D_X$, $\hat{\mu}^{1}=\mu^{1}$, $\hat{\boldsymbol{\lambda}}_{1}=\boldsymbol{\lambda}_{1}, J=1$ and $ k=1$.

  \WHILE {$k \leq \text{maxiter}$}
  \STATE Aggregator solves $\mathcal{D}_{0} (\boldsymbol{\hat{\lambda}}_{k})$ and broadcasts $\boldsymbol{\hat{\lambda}}_{k}$ and $\hat{\mu}^{k}$ to the households which solve $\mathcal{D}_{i,\hat{\mu}^{k}} (\boldsymbol{\hat{\lambda}}_{k})$ and return $\boldsymbol{x}_{\boldsymbol{i},\hat{\mu}^{k},\boldsymbol{\hat{\lambda}}_{k}}$ to the aggregator.
	\STATE Aggregator computes $\nabla \mathcal{D}_{\hat{\mu}^{k},\kappa^{k}}(\hat{\boldsymbol{\lambda}}_{k})$ and the primal $\mathcal{P}_{r}^{k}=\sum_{ t\in \mathcal{T}} C^{t}(\sum_{ i\in \mathcal{I} \setminus 0} x^{t}_{\boldsymbol{i},\hat{\mu}^{k},\boldsymbol{\hat{\lambda}}_{k}})$.
  \STATE Aggregator computes $L_{\mu^{k},\kappa^{k}}^{k}=\frac{\left\|A_c\right\|^2}{\mu^{k}}+\kappa^{k}$ and updates $\boldsymbol{\lambda}_{k+1}$ and $\hat{\boldsymbol{\lambda}}_{k+1}$.
	\STATE Aggregator updates $\mu^{k+1}$ and $\kappa^{k+1}$ and sets $\{\hat{\mu}^{k+1}=\mu^{k+1}: \hat{\mu}^{k+1} \geq \hat{\mu}^{\text{min}}\}$.
	\STATE $k \leftarrow k + 1 $.
	\ENDWHILE
 \renewcommand{\algorithmicensure}{\textbf{Output:}}
 \ENSURE  Aggregator finds the best recovered primal solution $\mathcal{P}_{r}^{J}$ along with $\boldsymbol{\hat{\lambda}}_{J}$, $\hat{\mu}^{J}$ and $\boldsymbol{x}_{\hat{\mu}^{J},\boldsymbol{\hat{\lambda}}_{J}}$ such that $J:=\{k:\mathcal{P}_{r}^{J}=\min \{ \{\mathcal{P}_{r}^{k}\}_{k \in \{1,\ldots,\text{maxiter}\} } \} \}$.
 \end{algorithmic} 
\label{algorithm}
 \end{algorithm}

\section{Numerical evaluation}\label{sec:results}

\begin{figure}[t]
\centering{
\psfrag{Iterations (k)}{\footnotesize Iterations (k) \normalsize}
\psfrag{Primal}{\footnotesize Primal \normalsize}
\psfrag{Dual}{\footnotesize Dual \normalsize}
\psfrag{Primal and Dual values}{\footnotesize $\mathcal{P}_{r}^{k}$ and $\mathcal{D}_{\hat{\mu}^{k},\kappa^{k}}(\hat{\boldsymbol{\lambda}}_{k})$ \normalsize}
\includegraphics[width=95mm] {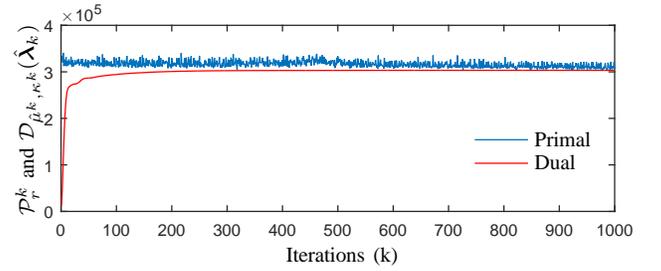}}
\caption{Evolution of the primal and dual objectives for the $I=1280$ case.}
\label{primaldual}
\end{figure}

Algorithm~\ref{algorithm} is tested on three scenarios, each with one aggregator and $640$, $1280$ and $2560$ households respectively, each with $10$ appliances scheduled over $T=24 \text{h}$ (as in \cite{FaithfulMDinDR}). In all three scenarios, Algorithm~\ref{algorithm} is initialized with $\boldsymbol{\lambda}_{1}=\boldsymbol{0}$, $\kappa^{1}=10$, $\text{maxiter}=1000$, $\alpha^{1}=3 \times 10^{-4}\left\|A_{c}\right\|^2 D_{X}$ and $\alpha^{\text{maxiter}}=8 \times 10^{-8}\left\|A_{c}\right\|^2 D_{X}$. The simulation results along with the corresponding parameter values are listed in Table~\ref{minimumcosts}. Table~\ref{minimumcosts} shows that the difference between the recovered best integer feasible solution $\mathcal{P}_{r}^{J}$ and the optimum solution $\mathcal{P}^{*}$ does not exceed $0.42 \%$, which corroborates the claim that a near-optimal solution can be recovered in a limited number of iterations. Finally, the evolution of $\mathcal{P}_{r}^{k}$ and $\mathcal{D}_{\hat{\mu}^{k},\kappa^{k}}(\hat{\boldsymbol{\lambda}}_{k})$ in the $I=1280$ case is displayed in Figure~\ref{primaldual}, which also shows a quick and smooth convergence of the dual objective $\mathcal{D}_{\hat{\mu}^{k},\kappa^{k}}(\hat{\boldsymbol{\lambda}}_{k})$ and a small duality gap upon termination.

\begin{table}[!t]
\centering
\caption{Difference between $\mathcal{P}_{r}^{J}$ and $\mathcal{P}^{*}$.}
\begin{tabular}{| r | c | c | c | c | c | c | c |}
\hline
$I$     & $\mathcal{P}_{r}^{J} (\$)$   & $\mathcal{P}^{*} (\$)$    &  Gap (\%)   	& $\hat{\mu}^{\text{min}}$  & $\kappa^{1}$     \\\hline
640  		& 76049.78  			& 75732.03     & 0.42    	   						& 0.0004 		 & 10   						\\\hline
1280  		& 304089.52  			& 302927.74    & 0.38    	  						& 0.001 		 & 10   						\\\hline
2560  		& 1216367.92  		  & 1211711.03    & 0.38  	 	 						& 0.0015 		 & 10   						\\\hline

\end{tabular}
\label{minimumcosts}
\end{table}

\section{Conclusion}\label{sec:conclusion}

The aim of this work is to implement a fast gradient algorithm applied to the double smoothed dual function of a DR problem comprising expressive household models and mixed-integer variables. This work also demonstrates how to recover a near-optimal solution in a fixed number of iterations and minimal parameter tweaking.

\bibliographystyle{IEEEtran}
{\footnotesize
\bibliography{ODMcitations}}

\begin{thebibliography}{1}
\providecommand{\url}[1]{#1}
\csname url@samestyle\endcsname
\providecommand{\newblock}{\relax}
\providecommand{\bibinfo}[2]{#2}
\providecommand{\BIBentrySTDinterwordspacing}{\spaceskip=0pt\relax}
\providecommand{\BIBentryALTinterwordstretchfactor}{4}
\providecommand{\BIBentryALTinterwordspacing}{\spaceskip=\fontdimen2\font plus
\BIBentryALTinterwordstretchfactor\fontdimen3\font minus
  \fontdimen4\font\relax}
\providecommand{\BIBforeignlanguage}[2]{{%
\expandafter\ifx\csname l@#1\endcsname\relax
\typeout{** WARNING: IEEEtran.bst: No hyphenation pattern has been}%
\typeout{** loaded for the language `#1'. Using the pattern for}%
\typeout{** the default language instead.}%
\else
\language=\csname l@#1\endcsname
\fi
#2}}
\providecommand{\BIBdecl}{\relax}
\BIBdecl

\bibitem{DRwithMIC}
S.-J. Kim and G.~Giannakis, ``Scalable and robust demand response with
  mixed-integer constraints,'' \emph{Smart Grid, IEEE Transactions on}, vol.~4,
  no.~4, pp. 2089--2099, Dec 2013.

\bibitem{FaithfulMDinDR}
S.~Mhanna, G.~Verbi\v{c}, and A.~Chapman, ``A faithful distributed mechanism
  for demand response aggregation,'' \emph{IEEE Trans. Smart Grid}, to be
  published.

\bibitem{subgradientmethods}
S.~Boyd, L.~Xiao, and A.~Mutapcic, ``Subgradient methods,'' \emph{lecture notes
  of EE392o, Stanford University, Autumn Quarter}, 2008.

\bibitem{intoncvxopt}
Y.~Nesterov, \emph{Introductory lectures on convex optimization}.\hskip 1em
  plus 0.5em minus 0.4em\relax Springer, 2004, vol.~87.

\bibitem{doublesmoothing}
O.~Devolder, F.~Glineur, and Y.~Nesterov, ``Double smoothing technique for
  large-scale linearly constrained convex optimization,'' \emph{SIAM Journal on
  Optimization}, vol.~22, no.~2, pp. 702--727, 2012.

\end{thebibliography}

\end{document}